\mag=\magstep1
\documentstyle{amsppt}

\topmatter
\title Index 1 covers of log terminal surface sigularities
\endtitle
\author Yujiro Kawamata
\endauthor

\rightheadtext{Index 1 covers}

\address Department of Mathematical Sciences, University of Tokyo, Komaba,
Meguro, Tokyo, 153-8914, Japan \endaddress
\email kawamata\@ms.u-tokyo.ac.jp\endemail


\subjclass 
14 E 22, 14 E 35, 14 J 17
\endsubjclass

\abstract
We shall investigate index 1 covers of $2$-dimensional 
log terminal singularities.
The main result is that the index 1 cover is canonical 
if the characteristic of the base field is different from $2$ or $3$.
We also give some counterexamples in the case of characteristic $2$ or $3$.
By using this result, we correct an error in [K2].
\endabstract

\endtopmatter

\document

\head 1. Introduction
\endhead

We fix an algebraically closed field $k$ and let $p$ be its
characteristic.
Let $S$ be a normal surface over $k$, let $P$ be a closed point of $S$, and
let $D$ be an effective and reduced Weil divisor on $S$ through $P$.
We consider the germ of the pair $(S, D)$ at $P$.
Let $\mu: S' \to S$ be an embedded resolution of the singularity for the 
pair $(S, D)$.
The {\it numerical pull-back} $\mu^*(K_S+D)$ is defined as a $\Bbb Q$-divisor 
on $S'$ such that the equality $(\mu^*(K_S+D)\cdot C) = ((K_S+D)\cdot \mu_*C)$ 
holds for any curve $C$ on $S'$. 
We write $\mu^*(K_S+D) = K_{S'} + D' + E$ for a $\Bbb Q$-divisor $E$ on $S'$,
where $D' = \mu_*^{-1}D$ is the strict transform of $D$.
The pair $(S, D)$ is said to be {\it log terminal} at $P$
if the coefficients of $E$ are strictly less than $1$.
It is called {\it canonical} if the coefficients of $E$ 
are non-positive.
$S$ is said to be canonical or log terminal if $(S, 0)$ is so
(cf. [KMM]).

The {\it index} $r$ of the pair $(S, D)$ 
is the smallest positive integer such that $r(K_S + D)$ is a 
Cartier divisor.
Let $\theta_0$ and $\theta$ be non-zero sections of 
$\Cal O_S(K_S + D)$ and $\Cal O_S(r(K_S + D))$, respectively.
Assume that $\theta$ generates $\Cal O_S(r(K_S+ D))$.
Let $L$ be the rational function field of $S$, 
and write $\theta = \alpha \theta_0^r$ for $\alpha \in L$.
The normalization $\pi: T \to S$ of $S$ in the field extension
$L(\alpha^{1/r})$ is called the {\it index 1 cover} of $S$ 
associated to the section $\theta$.
We note that the index 1 cover depends on the choice of $\theta$.
The index 1 cover is called the {\it log canonical cover} in [KMM], 
but its construction is not canonical at all, and 
in order to avoid a confusion, we use instead this terminology.

If $p$ does not divide the degree $r$ of the 
morphism $\pi$, then $\pi$ is etale over $S \setminus \{P\}$, and the 
pair $(T, D_T)$ for $D_T = \pi^*D$ is known to be caninical (cf. [KMM]).
But if $p$ divides $r$, then $\pi$ is inseparable, and the situation is
totally different.
The following is the main result of this paper.

\proclaim{Theorem 1}
Let $S$ be a normal surface over an algebraically closed field $k$
of characteristic $p \ne 2,3$, 
$D$ a reduced curve, and $P$ a closed point such that the pair $(S, D)$ is
log terminal of index $r$ at $P$.
Let $\theta$ be a nowhere vanishing section of $\Cal O_S(r(K_S + D))$, 
and $\pi: T \to S$ the index 1 cover associated to $\theta$.
Set $D_T = \pi^*D$ and $Q = \pi^{-1}(P)$.
If $\theta$ is chosen to be general enough, 
then $(T, D_T)$ is canonical at $Q$.
\endproclaim

We have counterexamples in the case of characteristic $2$ or $3$
(Example 5).
By using the above result, we shall correct an error in [K2] in \S 3.
We would like to thank Professor K. Matsuki for pointing out this error.

\head 2. Proof of Theorem 1
\endhead

We keep the notation of the introduction.

\proclaim{Lemma 2}
Assume that a pair $(S, D)$ is log terminal at a point $P$.
Then $(S, P)$ is a rational singularity.
\endproclaim

\demo{Proof}
Since $S$ is also log terminal, we may assume that $D = 0$.
Let $\mu: S' \to S$ be the minimal resolution and 
write $\mu^*K_S = K_{S'} + E$ for a $\Bbb Q$-divisor $E$.
Let $Z$ be the fundamental cycle.
Since the coefficients of $E$ are non-negative and less than $1$,
the divisor $Z - E$ is effective and its support is the whole exceptional 
locus of $\mu$.
Hence $(Z^2) + (Z \cdot K_{S'}) = (Z \cdot (Z - E)) < 0$.
\hfill $\square$
\enddemo

It follows that the index $r$ of the pair $(S, D)$ 
is equal to the smallest positive integer
such that $rE$ becomes a divisor.

Let us consider the divisor on $S'$
which is the sum of the exceptional locus of $\mu$ 
and the strict transform of $D$. 
Then we can classify the dual graphs of these divisors ([K1, TM]).
We note that this classification is purely numerical and characteristic free.
For example, the dual graphs for canonical singularities are Dynkin diagrams 
of type $A$, $D$ or $E$. 
The dual graphs of log terminal singularities are the same as those of
quotient singularities in characteristic $0$, but they are 
not necessarily quotient singularities in general.

We assume that $p \vert r$ from now on.
We start the proof of Theorem 1
with the calculation of the log canonical divisor on an index 
$r/p$ cover of a log termnal singularity of index $r$.
For a $1$-form $\omega$ on a normal variety $S$, we denote by
$\text{div}_S(\omega)$ the divisorial part of its zero or pole.

\proclaim{Lemma 3}
Let $S$ be a normal affine surface, 
$D$ a reduced curve, $P$ a closed point such that the pair $(S, D)$ is
log terminal of index $r$ at $P$.
Assume that $p \ne 2,3$ and $p$ divides $r$.
Let $\theta_0$ and $\theta$ be non-zero sections of 
$\Cal O_S(K_S + D)$ and $\Cal O_S(r(K_S + D))$, respectively.
Assume that $\theta$ is nowhere vanishing, and write
$\theta = \alpha \theta_0^r$ for $\alpha \in L$, the rational function field 
of $S$.
Assume that $\text{div}_S(d\alpha) = \text{div}_S(\alpha)$ as Weil divisors on
$S$.
Let $\pi: \widetilde S \to S$ be the normalization of $S$ in the field 
$\widetilde L = L(\alpha^{1/p})$. 
Set $\widetilde D = \pi^*D$ and $\widetilde P = \pi^{-1}(P)$.
Then $(\widetilde S, \widetilde D)$ is again log terminal at $\widetilde P$
of index $r/p$.
Moreover, if $p^2 \vert r$, then 
$\text{div}_{\widetilde S}(d\alpha^{1/p}) = 
\text{div}_{\widetilde S}(\alpha^{1/p})$ 
as Weil divisors on $\widetilde S$.
\endproclaim

\demo{Proof}
Since $\Cal O_S(\frac rp(K_S + D))$ is not invertible, 
we have $\alpha \not\in L^p$, and 
$\widetilde L/L$ is a purely inseparable extension of degree $p$.
By [RS, Proposition 2], we have
$$
K_{\widetilde S} = \pi^*(K_{S} - (1 - 1/p)\text{div}_S(d\alpha)).
$$
Since $\text{div}_S(d\alpha) = \text{div}_S(\alpha) \sim 0$, 
we have $K_{\widetilde S} = \pi^*K_{S}$ by a different choice of 
the identification.
Therefore, $\widetilde \theta = \alpha^{1/p}\pi^*\theta_0^{r/p}$
is a nowhere vanishing section of 
$\Cal O_{\widetilde S}(\frac rp (K_{\widetilde S} + \widetilde D))$, and
the index of the pair $(\widetilde S, \widetilde D)$ is $r/p$.

Let $\mu: S' \to S$ be a projective birational morphism from a smooth
surface, and
$D' = \mu_*^{-1}D$ the strict transform of $D$.
Since $(S, D)$ is log terminal, we can write
$$
\mu^*(K_S + D) = K_{S'} + D' + \sum_j a_j C_j
$$
with $a_j < 1$ for each irreducible component $C_j$ of the exceptional 
locus $C$ of $\mu$.
By the adjunction, the $C_j$ are isomorphic to $\Bbb P^1$, and intersect
transversally.

Let $\pi': \widetilde S' \to S'$ be the normalization in $\widetilde L$, 
$\widetilde \mu: \widetilde S' \to \widetilde S$ the induced birational 
morphism, and $\widetilde D' = \widetilde \mu_*^{-1}\widetilde D
= \pi^{\prime *}D'$.
We can write 
$$
\widetilde \mu^*(K_{\widetilde S} + \widetilde D) 
= K_{\widetilde S'} + \widetilde D' 
+ \sum_j \widetilde a_j \widetilde C_j,
$$
where the $\widetilde C_j$ are prime divisors such that 
$\pi'(\widetilde C_j) = C_j$. 
We know that $\pi^{\prime *}C_0 = \widetilde C_0$ or $p\widetilde C_0$.
We shall prove that $\widetilde a_j < 1$ for all $j$ and for any $\mu$.

By [RS, Proposition 2] again, we have
$$
K_{\widetilde S'} = \pi^{\prime *}(K_{S'} 
- (1 - 1/p)\text{div}_{S'}(d\alpha)).
$$
Since $\mu^*\text{div}_S(\alpha) = \text{div}_{S'}(\alpha)$, we have
$$
\sum_j \widetilde a_j \widetilde C_j
= \pi^{\prime *}(\sum_j a_jC_j - (1 - 1/p)(\text{div}_{S'}(\alpha)
- \text{div}_{S'}(d\alpha))).
$$

Let $G = \text{div}_S(\theta_0)$, and $G' = \mu_*^{-1}G$ its strict transform.
Then we can write 
$$
\align
\text{div}_{S'}(\alpha) + rG' &= \sum_j m_jC_j \\
\text{div}_{S'}(d\alpha) + rG' &= \sum_j m'_jC_j
\endalign
$$ 
for some $m_j, m_j' \in \Bbb Z$.  
Thus
$$
\widetilde a_j \widetilde C_j
= (a_j - (1 - 1/p)(m_j - m'_j))\pi^{\prime *}C_j.
$$
Since $G \sim K_S + D$, there exists a divisor $F$ supported on $C$ 
such that $F + G' \sim K_{S'} + D' $.
Then $\text{div}_{S'}(\alpha) + rF + rG' + \sum_j a_jrC_j$ is 
numerically trivial, hence $m_j + a_jr \equiv 0 \,(\text{mod }r)$.

Let us fix an irreducible component of $C$, say $C_0$.
We consider $2$ cases (we shall prove later that these are the only cases
provided that $p \ne 2,3$):

{\it Case 1.} 
We assume that $p$ does not divide $m_0$.  

We take a general closed point $P'$ on $C_0$.  Let $(x, y)$ be local 
coordinates such that $C_0 = \text{div}(x)$ near $P'$.
We can write $\alpha = ux^{m_0}$ near $P'$ such that $u(P') \ne 0$.
Then we have $d\alpha = x^{m_0-1}(m_0udx + xdu)$, hence $m'_0 = m_0 - 1$.
Since $p$ does not divide $m_0$, there are integers $s, t$ such that
$ps + m_0t = 1$.  
Then $\text{div}_{S'}(x^{ps}\alpha^t) = C_0$ near $P'$, hence
$\pi^{\prime *}C_0 = p\widetilde C_0$ with
$\text{div}_{\widetilde S'}(x^s\alpha^{t/p}) = \widetilde C_0$
near $\widetilde P' = \pi^{\prime -1}(P')$.
Therefore, $\widetilde a_0 = p(a_0 - 1 + 1/p) < 1$. 

{\it Case 2.}
We assume that $p \vert m_0$. 
In this case, we assume in addition that $C_0$ intersects at most $2$ other 
components of $C$,
say $C_1$ and $C_2$ ($C_2$ may not exist).
Moreover, we assume that $p$ does not divide $m_1$.  

Since $m_1 + m_2 \equiv ra_1 + ra_2 \equiv 0 \, (\text{mod }p)$,
$C_2$ necessarily exists and $p$ does not divide $m_2$. 
Let $P'$ be an arbitrary closed point on $C_0$ except $P'_i = C_0 \cap C_i$ for
$i = 1, 2$, and $(x, y)$ local coordinates such that $C_0 = \text{div}(x)$ 
near $P'$. 
We can write $\alpha = uv^px^{m_0}$ near $P'$ in such a way that 
$u_0 = u\vert_{C_0}$ is a rational function on $C_0$ such that
$\text{div}_{C_0}(u_0) = m_1P'_1 + m_2P'_2 + pQ'$ for some divisor $Q'$ on 
$C_0$ whose support does not contain $P'$.
Since the $m_i$ are not divisible by $p$, we have $du_0 \ne 0$.
Thus $\text{deg}(du_0) = -2$, and we have
$\text{div}_{C_0}(du_0) = (m_1-1)P'_1 + (m_2-1)P'_2 + pQ'$.
Therefore, $du_0$ does not vanish at $P'$.
Since $d\alpha = v^px^{m_0}du$ near $P'$, we have $m'_0 = m_0$.
Moreover, $u_0 - u_0(P')$ gives a local coordinate of $C_0$ at $P'$.
Hence $(\pi^{\prime *}x, \pi^{\prime *}(u - u(P'))^{1/p})$ give local
coordinates at $\widetilde P' = \pi^{\prime -1}(P')$.
Thus $\pi^{\prime *}C_0$ is reduced, and $\widetilde S'$ is smooth at 
$\widetilde P'$.
In particular, $\pi^{\prime *}C_0 = \widetilde C_0$ and 
$\widetilde a_0 = a_0 < 1$. 

We shall prove that any irreducible component $C_0$ of $C$ satisfies 
the assumptions of one of the above two cases. 
First, we consider the case in which $\mu = \mu_0: S' = S'_0 \to S$ 
coincides with the minimal resolution. 

Assuming that $C_0$ intersects $3$ other components, say $C_1$, $C_2$, $C_3$,
we shall prove that we have Case 1 for $C_0$.
Assume the contrary that $p \vert m_0$. Then $p \vert a_0r$.
In the case in which the dual graph for $S'$ is of type $D$, 
we have $(C_1^2) = (C_2^2) = -2$ after the permutation of the indices.
Then we calculate that $a_1 = a_2 = a_0/2$.
Since $m_j + a_jr \equiv 0 \,(\text{mod }r)$ and $p \ne 2$, 
we have $p \vert m_1$ and $p \vert m_2$.
Since $p \vert m_1 + m_2 + m_3 - (C_0^2)m_0$, we have $p \vert m_3$.
Then we have $p \vert m_4$ for an irreducible component $C_4$ which intersects
$C_3$.  In this way, we conclude that $p \vert m_j$ for all $j$.
It follows that $\frac rp(K_{S'} + D' + \sum_ja_jC_j)$ is a divisor on $S'$.
Since this divisor is numerically trivial and $S$ is a rational singularity, 
it is a pull back of a divisor on $S$, 
a contradiction with the assumption that $r$ is the index.

In the case in which the dual graph for $S'$ is of type $E$, 
we have two cases after the permutation of the 
indices: (i) $(C_1^2) = (C_2^2) = - 2$ and $C_2$ intersects another 
irreducible component $C_4$ such that $(C_4^2) = -2$ while $C_1$ does
not intersects other components, 
or (ii) $(C_1^2) = -2$, $(C_2^2) = - 3$, and $C_1$ and $C_2$ intersect no
other irreducible components. 
We have $a_1 = a_0/2$ and $a_2 = 2a_0/3$ in the former case, 
and $a_1 = a_0/2$ and $a_2 = (a_0 + 1)/3$ in the latter.
Since $p \ne 2, 3$, we have $p \vert m_1$ and $p \vert m_2$, 
and obtain a contradiction as before.

If we assume that $p \vert m_0$,
then by the above argument, 
$C_0$ intersects at most $2$ other components of $C$,
say $C_1$ and $C_2$ ($C_2$ may not exist).
Suppose that $p \vert m_1$.  
Then $C_1$ intersects at most $1$ other component, say $C_3$, and that
$p \vert m_3$. 
Moreover, if $C_2$ exists, then we have also $p \vert m_2$.  
Then we have $p \vert m_j$ for all $j$ as before, a contradiction. 
Therefore, we have Case 2 for $C_0$.

Next, we consider the general case.
Let $\mu$ and $j$ be arbitrary.
If the center of $C_j$ on the minimal resolution $S'_0$ is a curve, then
the above argument showed our assertion.
Assume that the center $Q$ of $C_j$ on $S'_0$ is a point. 
We have 3 cases: (a) $Q$ is contained in only one irreducible component $C_0$ 
of $C$ such that $p \vert m_0$, 
(b) $Q$ is contained in only one irreducible component $C_0$ 
of $C$ such that $p \not\vert m_0$, 
(c) $Q$ is contained in two irreducible components $C_0$ and $C_1$ of 
$C$ such that $p \not\vert m_0$.

In the case (a), the covering $\widetilde S'_0$ is smooth 
at the point $\widetilde Q$ above $Q$ by the argument of Case 2,
hence we obtain $\widetilde a_j < 1$ after any sequence of blow-ups of
$\widetilde S'_0$ above $\widetilde Q$.

In the case (b) or (c), we replace $S'_0$ by its blow-up at $Q$, and we 
obtain again (a), (b) or (c), where $C_0$ may be equal to $C_j$ in the 
cases (a) or (b). In the latter cases, we have Case 1 or 2 and are done.
Therefore, after finitely many blow-ups, we deduce that $\widetilde a_j < 1$.
We note that the above proof also showed that $\widetilde D$ is reduced.

Finally, in order to prove the last statement, we claim that 
$$
p\text{div}_{\widetilde S}(d\alpha^{1/p}) = 
\pi^*(\text{div}_S(d\alpha)).
$$
We shall check this equality at all but finitely many points on $S$.
As in the proof of [RS, Proposition 2], we may assume that 
there exist local coordinates $(x, y)$ of the completion of $S$ 
at the point $Q$ such that 
$(\widetilde x, \widetilde y)$ with $\widetilde x = \pi^*x$ and
$\widetilde y = \pi^*y^{1/p}$ give local
coordinates of the completion of $\widetilde S$ at 
$\widetilde Q = \pi^{-1}(Q)$.
We can write $\alpha = u^r\sum_{i=0}^{p-1}c_i^py^i$ for $u \in L$ and 
$c_i \in \hat L$, where $\hat L$ is the fraction field of the completed local 
ring, such that $\text{div}_S(\alpha) = r\text{div}_S(u)$.
Since $\text{div}_S(d\alpha) = \text{div}_S(\alpha)$, we may assume that
$c_1(Q) \ne 0$.
Since 
$$
d\alpha^{1/p} = u^{r/p}\sum_{i=0}^{p-1}
(ic_i\widetilde y^{i-1}d\widetilde y + \widetilde y^idc_i),
$$
we have $p\text{div}_{\widetilde S}(d\alpha^{1/p}) 
= r\text{div}_{\widetilde S}(u)$.
\hfill $\square$
\enddemo

\demo{Proof of Theorem 1}
Since $\theta$ is chosen to be general enough,
we deduce that $\text{div}(d\alpha) = \text{div}(\alpha)$ 
if we replace $S$ by a suitable neighborhood of $P$
by the dimension count argument as in p. 472 of [K2].
We apply Lemma 3 until the index becomes coprime to $p$, then 
apply the usual argument to obtain our assertion (cf. [KMM]).
\hfill $\square$
\enddemo

\definition{Remark 4}
(1) The formula for $K_{\widetilde S}$ depends on the choice of 
$\alpha^{1/p}$ which generates the field extension
$\widetilde L/L$.
This choice is equivalent to the splitting of a free $L$-module $\widetilde L$ 
as 
$$
\widetilde L = \bigoplus_{m=0}^{p-1} L\alpha^{m/p}.
$$ 
The construction of index 1 cover as in [K2] uses this kind of splitting 
explicitly and thus there is a canonical divisor formula.

(2) Lemma 3 is still true in the case of characteristic $2$ or $3$ if 
the minimal resolution diagram of $S$ is of type $A$.
\enddefinition

\definition{Example 5}
(1) Let $\mu: S' \to S$ be the minimal resolution of a log terminal 
singularity over an algebraically closed field of characteristic $p = 2$.
Assume that the dual graph of the exceptional divisors $C$ is of type $D$ as 
follows: $C = C_1 + C_2 + C_3 + C_4$ with 
$(C_1^2) = (C_2^2) = (C_3^2) = -2$, $(C_4^2) = - 3$,
$(C_1 \cdot C_2) = (C_1 \cdot C_3) = (C_1 \cdot C_4) = 1$, and 
other intersection numbers are $0$.
We note that a surface $S$ as above can be constructed by
blowing up suitably a smooth surface and then contracting some of the 
exceptional divisors.
We have
$$
\mu^*K_S  = K_{S'} + \frac 12 C_1 + \frac 14 C_2 + \frac 14 C_3 + \frac 12 C_4. 
$$
$S$ is a rational triple point and the index $r = 4$.
Let $\pi: T @>{\pi_2}>> \widetilde S @>{\pi_1}>> S$ 
be the index 1 cover associated to a general section $\theta$ of 
$\Cal O_S(4K_S)$, where $\pi_1$ and $\pi_2$ are purely inseparable morphisms 
of degree $2$.
We claim that $T$ is not log terminal.

Indeed, as in Case 2 in the proof of Lemma 3,
since $2 \vert m_1$,
we can write $\alpha = uv^2x^{m_4}$ near a general closed point $P'$ 
of $C_4$ in such a way that $u_4 = u\vert_{C_4}$
is a rational function on $C_4$ 
such that $\text{div}_{C_4}(u_4) = 2Q'$ for some divisor $Q'$ on $C_4$.
Thus $u_4 = v_4^2$ for some rational function $v_4$ on $C_4$.
It follows that the natural morphism
$\widetilde C_4 \to C_4$ is birational, hence 
$\pi^{\prime *}C_4 = 2\widetilde C_4$.
Since $m'_4 \ge m_4$, we have $\widetilde a_4 \ge 1$. 
If we denote by $b_j$ the coefficients for $K_T$ in a suitable way,
then we deduce that $b_4 \ge 1$ by the same argument as in the proof 
of Lemma 3.

(2) Let $\mu: S' \to S$ be the minimal resolution of a log terminal 
singularity over an algebraically closed field of characteristic $p = 3$.
Assume that the dual graph of the exceptional divisors $C$ is of type $E_6$ as 
follows: $C = C_1 + C_2 + C_3 + C_4 + C_5$ with 
$(C_1^2) = (C_2^2) = (C_3^2) = (C_4^2) = -2$, $(C_5^2) = - 3$,
$(C_1 \cdot C_2) = (C_1 \cdot C_3) = (C_3 \cdot C_4) = (C_1 \cdot C_5) = 1$, 
and other intersection numbers are $0$.
Then we have
$$
\mu^*K_S = K_{S'} + \frac 23 C_1 + \frac 13 C_2 + \frac 49 C_3 + 
\frac 29 C_4 + \frac 59 C_5. 
$$
$S$ is a rational quintuple point and $r = 9$.
Let $\pi: T \to \widetilde S \to S$ be the 
index 1 cover associated to a general section $\theta$ of $\Cal O_S(9K_S)$. 
We claim that $T$ is not log terminal.
Indeed, we have $b_2 \ge 1$ as in (1).
\enddefinition

\head 3. Correction to [K2]
\endhead

Kenji Matsuki pointed out that the proof of Theorem 3.1 of [K2] 
is insufficient 
because the calculation in the middle of p.473 is wrong.
We shall replace the proof of Theorems 3.1 and 4.1 
of [K2] by a different argument and prove them 
under the additional assumption that the 
residue characteristic is different from $2$ or $3$.
We note that it is still an open question in the case of characteristic 
$2$ or $3$.

\demo{Proof of Theorems 3.1 and 4.1 of [K2] in the case $p \ne 2,3$}
We prove Theorem 4.1 by a slightly modified argument.
Theorem 3.1 follows a posteriori from Theorem 4.1.
We use the notation in Theorem 3.1;
let $f: X \to \Delta = \text{Spec }A$ be a family satisfying Assumption 1.1.
Let $p$ be the characteristic of the residue field at the closed point of
$\Delta$. We assume that $p \ne 2,3$.
We take a closed point $P \in X$ of index $r$.
The index 1 cover $\pi: Y \to X$ is constructed by using a general 
section $\theta$ of $\Cal O_X(rK_{X/\Delta})$ as 
$$
\pi_*\Cal O_Y \cong \bigoplus_{m=0}^{r-1} \Cal O_X(-mK_{X/\Delta})t^m, \,
t^r = \theta.
$$
We shall prove that the closed fiber $Y_s$ is canonical or 
normal crossing, but we do not prove that the singularity of $Y$ is isolated
at this point.

First, assume that the closed fiber $X_s$ of $X$ is irreducible. 
Since $\Cal O_X(-mK_{X/\Delta}) \otimes \Cal O_{X_s}
\cong \Cal O_{X_s}(-mK_{X_s})$ by Assumption 1.1 (6), 
$Y_s$ is isomorphic to the index 1 cover of $X_s$
constructed by using the restriction of $\theta$ to $X_s$. 
By Theorem 1, $Y_s$ is canonical. We can prove that
its completed local ring at $Q$ is
isomorphic to the completion of $A[x_1,x_2,x_3]/(F)$ with 
$\text{ord}(F_s) \le 2$ as in the original proof of Theorem 4.1, where 
the action of $\mu_r$ on the coordinates $(x_1,x_2,x_3)$ is given by
$x_i \mapsto \zeta^{a_i} \otimes x_i$ ($i = 1,2,3$).
Since $\Cal O_X(-K_{X/\Delta})$ is not invertible, there exists at least $2$
coordinates whose weights $a_i$ are coprime to $r$.
Let $x_1, \ldots, x_c$ ($c = 2$ or $3$) be such coordinates. 
Since $\theta = t^r$ never vanishes and the natural homomorphism
$\Cal O_X(-K_{X/\Delta})^{\otimes r} \to \Cal O_X(-rK_{X/\Delta})$ is 
surjective outside $\{P\}$, 
we have $\{x_1 = \cdots = x_c = F = 0\} = \{Q\}$.
It follows that all the $a_i$ are coprime to $r$ and $F$ contains a term 
in $A$.
Thus $F_s$ is $\mu_r$-invariant, and $\text{ord}(F_s) = 2$.
If $F_s$ contains a term of the form $x_1x_2$, then we are done.
If it contains $x_1^2$ and there are no other terms of degree $2$, then
$r = 2$. But there is a term of degree $3$ in $F_s$, a contradiction.

Next, assume that $X_s$ is reducible.
Let $X_{s,i}$ ($1 \le i \le d$) be its irreducible components.
Since the $X_{s,i}$ are $\Bbb Q$-Cartier divisors and the pairs 
$(X_{s,i}, D_i)$ for $D_i = \sum_{j\ne i}X_{s,j}\cap X_{s,i}$ 
are log terminal, we have $d = 2$ or $3$.
Let $r_i$ be the indices of the $(X_{s,i}, D_i)$.
If $d = 3$, then $r_i = 1$ for all $i$, and there are nowhere vanishing 
sections $\theta_i$ of $\Cal O_{X_{s,i}}(K_{X_{s,i}} + D_i)$ which coincide 
each other on the double locus of $X_s$ to give a nowhere vanishing section 
$\theta_{X_s}$ of $\Cal O_{X_s}(K_{X_s})$.
Here we used the assumption that $\Cal O_{X_s}(K_{X_s})$ has depth $2$ at $P$.
Therefore, $r = 1$, a contradiction.

We consider the case $d = 2$.
$\theta$ induces a section $\theta_{X_s}$ of 
$\Cal O_{X_s}(rK_{X_s})$ and the sections $\theta_i$ of the 
$\Cal O_{X_{s,i}}(r(K_{X_{s,i}}+D_i))$. 
Thus $r_i \vert r$.
We write $r = r'p^f$ and $r_i = r'_ip^{f_i}$ with $(r', p) = 1$ and $(r'_i, p) 
= 1$.
We can construct a covering $\pi': X' \to X$ of degree $r'$ by
$$
\pi'_*\Cal O_{X'} \cong \bigoplus_{m=0}^{r'-1} \Cal O_X(-mp^fK_{X/\Delta})t^m,
\, t^{r'} = \theta.
$$
Then $X'_{s,i} = \pi^{\prime -1}X_{s,i}$ is a union of $r'/r'_i$ 
prime divisors which intersect only at a point $\pi^{\prime -1}(P)$.
Since $X'_{s,i}$ supports a Cartier divisor on $X'$, it follows that 
$r' = r'_i$ for $i = 1,2$.

Let $\theta_0$ be a section of $\Cal O_{X_s}(K_{X_s})$ which does not vanish
identically along the double locus $D$ of $X_s$.
We write $\theta_{X_s} = \alpha \theta_0^r$ as in Lemma 3.
Since $\theta$ is general, we may assume that $\alpha_D = \alpha \vert_D 
\not\in L_D^p$,
where $L_D$ is the rational functin field of $D$.
Let $Y_{s,i} = \pi^{-1}(X_{s,i})$.
We can extend Lemma 3 and apply it to the induced morphism 
$\pi_i: Y_{s,i} \to X_{s,i}$ 
even if $r_i$ might be smaller than $r$, 
because $L_D(\alpha_D^{1/p^f})/L_D$ is a purely inseparable
field extension.
Since $\Cal O_X(-mK_{X/\Delta}) \otimes \Cal O_{X_{s,i}} 
\cong \Cal O_X(-m(K_{X_{s,i}}+D_i))$ on $X_{s,i} \setminus \{P\}$,
$Y_{s,i}$ is smooth possibly except at $Q$, and $\pi_i^*D_i$ is a reduced 
smooth divisor on $Y_{s,i} \setminus \{Q\}$.
Thus $Y_s \setminus \{Q\}$ is a normal crossing divisor on $Y \setminus \{Q\}$.
Since $Y_s$ has depth $2$, we conclude that
the completed local ring of $Y$ is
isomorphic to the completion of $A[x_1,x_2,x_3]/(F)$ with 
$F_s = x_1x_2$ as in the original proof of Theorem 4.1. We may assume that
the action of $\mu_r$ on the coordinates $(x_1,x_2,x_3)$ is given by
$x_i \mapsto \zeta^{a_i} \otimes x_i$ ($i = 1,2,3$), because the ideal 
$(F_s)$ is preserved by this action.
By the same reason as in the case where $X_s$ is irreducible,
all the $a_i$ are coprime to $r$ and $F$ contains a term in $A$ 
so that $F = x_1x_2 + \tau$.
Since the $X_{s,i}$ are Cartier divisors on $X \setminus \{P\}$, so are
the $Y_{s,i}$ on $Y \setminus \{Q\}$.
Therefore, $Y \setminus \{Q\}$ is regular, and 
$\tau$ is a generator of the maximal ideal of $A$.
\hfill $\square$
\enddemo

\enddocument